\documentclass[10pt]{article}
\usepackage{amsmath,amssymb,fullpage}
\usepackage[T1]{fontenc}

\def\qed{\hfill $\Box$}
\def\one #1{1_{\{#1\}}}

\def\diag{\mbox{diag}}
\def\diag #1{{\rm diag}\left(#1\right)}
\def\1{\mbox{\bf 1}}
\def\cadlag{c\`{a}dl\`{a}g\ }
\def\n{\nonumber \\}
\def\eq #1{(\ref{eq:#1})}

\newcommand{\proof}{\noindent {\bf Proof:\ }}

\newtheorem{Theorem}{Theorem}
\newtheorem{Lemma}{Lemma}

\newtheorem{Corollary}{Corollary}
\newtheorem{Remark}{Remark}

\newtheorem{Assumption}{Assumption}

\begin{document}
\title{Unifying the Dynkin and Lebesgue-Stieltjes formulae}

\author{
Offer Kella\thanks{Department of Statistics; The Hebrew
University of Jerusalem; Jerusalem 91905; Israel
({\tt Offer.Kella@huji.ac.il}).} \thanks{Supported by grant No. 1462/13 from the Israel Science Foundation and
the Vigevani Chair in Statistics.}
\and Marc Yor\thanks{Laboratoire de Probabilit\'es et Mod\`eles al\'eatoires; Universit\'e Pierre et Marie Curie; Bo\^ite courrier 188; 75252 Paris Cedex 05; France ({\tt deaproba@proba.jussieu.fr})}
}
\date{June 5, 2016}
\maketitle
\begin{abstract}
We establish a local martingale $M$ associate with $f(X,Y)$ under some restrictions on $f$, where $Y$ is a process of bounded variation (on compact intervals) and either $X$ is a jump diffusion (a special case being a L\'evy process) or $X$ is some general (c\'adl\'ag metric space valued) Markov process. In the latter case $f$ is restricted to the form $f(x,y)=\sum_{k=1}^K\xi_k(x)\eta_k(y)$. This local martingale unifies both Dynkin's formula for Markov processes and the Lebesgue-Stieltjes integration (change of variable) formula for (right continuous) functions of bounded variation. For the jump diffusion case, when further relatively easily verifiable conditions are assumed then this local martingale becomes an $L^2$ martingale. Convergence of the product of this Martingale with some deterministic function (of time) to zero both in $L^2$ and a.s. is also considered and sufficient conditions for functions for which this happens are identified.
\end{abstract}


\bigskip
\noindent {\bf AMS Subject Classification (MSC2000):} Primary 60G44; Secondary 60J25, 60G51, 60K30.

\section{Introduction}
Suppose that $X=\{X_t|\ t\ge 0\}$ is a \cadlag Markov process taking values in some metric space where the notion of \cadlag is well defined with respect to some filtration satisfying the usual conditions and that $Y=\{Y_t|\ t\ge0\}$ is a $\mathbb{R}^r$ valued c\`adl\`ag, adapted process of finite variation on finite intervals (FV). Assume that $X$ has an extended generator $\mathcal{A}$ such that for any continuous $\xi$ in its domain (so that $\xi(X_t)$ is also \cadlag and adapted) we have that
\begin{align}
\xi(X_t)=\xi(X_0)+\int_0^t\mathcal{A}\xi(X_s)ds+M^\xi_t
\end{align}
where $M^\xi$ is a local Martingale (Dynkin's formula). If $\eta$ is a continuously differentiable function, then the Lebesgue-Stieltjes integration formula reads
\begin{align}\eta(Y_t)=\eta(Y_0)+\int_0^t\nabla\eta(Y_s)^TdY^c_s+\sum_{0<s\le t}\Delta\eta(Y_s)\end{align}
where $\Delta Y_s=Y_s-Y_{s-}$, $Y^c_t=Y_t-\sum_{0<s\le t}\Delta Y_s$ and $\nabla$ is the gradient operator.

As $\xi(X)$, $\eta(Y)$ are (c\`adl\`ag) semimartingales, integration by parts (\cite{p2004}, p.~68, Cor.~2) gives
\begin{align}
\xi(X_t)\eta(Y_t)&=\xi(X_0)\eta(Y_0)+\underbrace{\int_{(0,t]}\eta(Y_{s-})d\xi(X_s)}_1+\underbrace{\int_{(0,t]}\xi(X_{s-})d\eta(Y_s)}_2+\underbrace{[\xi(X),\eta(Y)]_t}_3\n
&=\xi(X_0)\eta(Y_0)+\underbrace{\int_0^t\mathcal{A}\xi(X_s)\eta(Y_s)ds+\int_{(0,t]}\eta(Y_{s-})dM^\xi_s}_1\\
&\ \ +\underbrace{\int_0^t\xi(X_s)\nabla\eta(Y_s)^TdY^c_s+\sum_{0<s\le t}\xi(X_{s-})\Delta f(Y_s)}_2+\underbrace{\sum_{0<s\le t}\Delta\xi(X_s)\Delta\eta(Y_s)}_3\n
&=\xi(X_0)\eta(Y_0)+\int_0^t\mathcal{A}\xi(X_s)\eta(Y_s)ds+\int_0^t\xi(X_s)\nabla\eta(Y_s)^TdY^c_s\n
&\ \ +\sum_{0<s\le t}(\xi(X_s)\eta(Y_s)-\xi(X_s)\eta(Y_{s-}))+\int_{(0,t]}\eta(Y_{s-})dM^\xi_s\ .\nonumber
\end{align}
 Letting $f(x,y)=\xi(x)\eta(y)$, note that for each fixed $y$, $f_y(x)=f(x,y)$ is in the domain of $\mathcal{A}$ and thus we will abuse the notation of $\mathcal{A}$ to denote an operator satisfying for each fixed $y$ that $\mathcal{A}f(x,y)=\mathcal{A}f_y(x)$. Also, for each fixed $x$, we denote $\nabla_yf(x,y)$ the gradient of $f_x(y)=f(x,y)$ with respect to $y$. With these conventions, noting that $\int_{(0,t]}\eta(Y_{s-})dM_s^\xi$ is a local martingale, we obtain
 \begin{Theorem}
For c\`adl\`ag adapted $X$, $Y$, where $X$ is Markov on some metric space with extended generator $\mathcal{A}$, $Y$ is $\mathbb{R}^r$ valued FV and $f(x,y)=\xi(x)\eta(y)$ where $\xi$ is continuous, in the domain of $\mathcal{A}$ and $\eta$ is continuously differentiable, then with the definition of $\mathcal{A}f$ and $\nabla_yf$ in the preceding paragraph we have that
\begin{align}\label{eq:M_t}
M_t&=f(X_t,Y_t)-f(X_0,Y_0)-\int_0^t\mathcal{A}f(X_s,Y_s)ds\\&\ \ -\int_0^t\nabla_yf(X_s,Y_s)^TdY^c_s-\sum_{0<s\le t}(f(X_s,Y_s)-f(X_s,Y_{s-}))\nonumber
\end{align}
is a local martingale.
\end{Theorem}

We note that if $X$ is a finite dimensional semimartingale, the sum part in the generalized (to discontinuous semimartingales) It\^o's formula is different from the simple form
\begin{align}
\sum_{0<s\le t}(f(X_s,Y_s)-f(X_s,Y_{s-}))
\end{align}
where we emphasize that $X_s$ (not $X_{s-}$) appears in both $f(X_s,Y_s)$ and $f(X_s,Y_{s-})$.

Assuming that \eq{M_t} is a local martingale for some continuous function $f$ for which $f(\cdot,y)$ is in the domain of $\mathcal{A}$ and $f(x,\cdot)$ is continuously differentiable, then when $f$ does not depend on $Y$ or $Y$ is constant, then the second line of \eq{M_t} is zero and \eq{M_t} reduces to Dynkin's formula. When $f$ does not depend on $X$ or $X$ is constant, then $\int_0^t\mathcal{A}f(X_s,Y_s)ds=M_t=0$ and then \eq{M_t} reduces to the Lebesgue-Stieltjes integration formula. When $X$ is Brownian motion and $Y$ is continuous then \eq{M_t} becomes It\^o's formula.

We immediately note that \eq{M_t} is a local martingale also for functions of the form $\sum_{k=1}^K\xi_k(x)\eta_k(y)$ where $\xi_k$ are continuous and in the domain of $\mathcal{A}$ and $\eta_k$ are continuously differentiable. Thus, it is plausible that with the aid of a Stone-Weierstrass type theorem the class of functions for which \eq{M_t} is a local martingale is far broader. As the main goal is to identify conditions under which \eq{M_t} is a martingale, we leave this approach to a future study. Therefore, in this paper we will show that \eq{M_t} is a local martingale for a relatively general special case where $X$ is a (possibly multivariate) jump-diffusion process, $f(\cdot,y)$ is twice differentiable, $f(x,\cdot)$ is differentiable and all (mixed) derivatives are continuous, also satisfying that the integral of $f(\cdot,y)$ with respect to some L\'evy kernel is finite for each relevant $y$. In addition, we will provide (for this case) relatively easily verifiable sufficient conditions for \eq{M_t} to be an $L^2$ martingale as well as satisfy $h(t)M_t\to0$ in $L^2$ or a.s. (for some $h$) under appropriate conditions (see Lemma~\ref{C}).

When $X$ is a real valued L\'evy process with triplet $(c,\sigma^2,\nu)$ and $Y$ is real valued, the operator $\mathcal{A}$ (for sufficiently nice functions is given by (e.g., Section~3.5.1 of \cite{applebaum})
\begin{align}
\mathcal{A}f(x,y)&=cf_x(x,y)+\frac{\sigma^2}{2}f_{xx}(x,y)\nonumber\\&\ \ +\int_{\mathbb{R}}(f(x+z,y)-f(x,y)-f_x(x,y)z\one{|z|\le 1})\nu(dz)
\end{align}
and in particular if we take $f_1(x,y)=\cos(\alpha(x+y))$ and $f_2(x,y)=\sin(\alpha(x+y))$, then for $f=f_1+if_2=e^{i\alpha(x+y)}$ we have that
\begin{align}\label{eq:Apsi}
\mathcal{A}f(x,y)=\psi(\alpha)e^{i\alpha(x+y)}
\end{align}
where
$\psi$ is the L\'evy exponent given by
\begin{align}\label{eq:psi}
\psi(\alpha)=ic\alpha-\frac{\sigma^2}{2}\alpha^2+\int_{\mathbb{R}}(e^{i\alpha z}-1-i\alpha z\one{|z|\le 1})\nu(dz)\ .
\end{align}
For this case $M$ is the (local) martingale from \cite{kw1992}. Originally, this process was shown to be only a local martingale, unless some further conditions were assumed, but it was discovered in \cite{kb2013} that it, as well as a generalized version of it, is in fact always an $L^2$ martingale and moreover $M_t/t\to0$ a.s. and in $L^2$. This latter result was needed to prove quite general decomposition results for on/off L\'evy storage processes and polling systems (see~\cite{bk2014}). In \cite{ak2000} a certain (but incomplete) generalization of \cite{kw1992} was established with various applications. The desire to place the results from \cite{kw1992,kb2013,ak2000} in a more general setting with formulas which seem easier to remember and straightforward to apply by experts and non-experts alike is what motivated the current study.

We note that the results from \cite{kw1992,ak2000} have been used in quite a few theoretical and applied studies (mainly in queueing, risk and finance) and extensively cited over the years. We find it unnecessary to give an exhaustive summary, as this can easily be found via a simple web search. Some book reference may be found in \cite{a2003,aa2010,w2002,k2006,k2013,lm2015,dm2015,kp2004}. A very small random sample of applications is, ({\em e.g.}): \cite{aap2004} (finance), \cite{f2005,psy2003} (risk), \cite{bps2001} (queueing), \cite{pr2002,ap2007} (theoretical).

The paper is organized as follows. In Section~\ref{Levy} we show in some detail how the main ideas work for the case where $X$ is a real valued L\'evy process and the process $Y$ is also one dimensional. In particular we show how the corresponding results from \cite{kw1992,kb2013} are obtained as immediate special cases. In Section~\ref{multijumpdiffusion} these ideas are generalized to multivariate jump diffusion $X$ and multivariate $Y$. In Section~\ref{markovadditive} we apply the results of Section~\ref{multijumpdiffusion} to certain Markov additive processes and generalize in more than one way corresponding results reported in \cite{ak2000}.

Since the important ideas can be found already in the L\'evy case, we thought that it would be easier to follow the derivation for this basic case and only then mention what is needed in order to generalize to the more general case. The added benefit is that it would make this more accessible to those who are only interested in the L\'evy case.

Throughout the paper, $\mathbb{R}_+$ is the set of nonnegative reals, $a\wedge b=\min(a,b)$, $a\vee b=\max(a,b)$ and a.s. abbreviates {\em almost surely}. Also, for semimartingale $U,V$, $[U,V]$ denotes the covariation process.

\section{L\'evy $X$, one dimensional $X$ and $Y$}\label{Levy}
Here we assume that $Y$ is FV and $X$ is a real valued L\'evy process with associated:
\begin{itemize}
\item L\'evy triplet $(c,\sigma^2,\nu(\cdot))$, where $c\in\mathbb{R}$, $\sigma\ge 0$,  $\int_{\mathbb{R}}(z^2\wedge 1)\nu(dz)<\infty$ and $\nu(\{0\})=0$.
\item Wiener process $W$.
\item Poisson random measure $N(dz,dt)$ with mean measure $\nu(dz)dt$.
\item $\tilde N(dz,dt)=N(dz,dt)-\nu(dz)dt$ (compensated Poisson measure).
\end{itemize}
That is,
\begin{align}\label{eq:LK}
X_t=ct+\sigma W_t+\int_{\mathbb{R}\setminus [-1,1]\times(0,t]}zN(dz,ds)+\int_{[-1,1]\times (0,t]}z\tilde N(dz,ds)
\end{align}
where it is well known that every real valued L\'evy process has such a decomposition (e.g., Thm.~42, p. 31 of \cite{p2004}).

If $f\in\mathcal{C}^{2,1}$ (e.g., \cite{ny2005}), as every L\'evy process and every FV process are semimartingales (recall that all processes are assumed adapted and c\`adl\`ag), according to the standard generalization of It\^o's lemma,
\begin{align}
f(X_{t},Y_{t})&=f(X_{0},Y_{0})+\int_{(0,t]}f_{x}(X_{s-},Y_{s-})dX_{s}+\int_{(0,t]}f_{y}(X_{s-},Y_{s-})dY_{s}\n
&\ +\frac{1}{2}\int_{0}^{t}f_{xx}(X_{s},Y_{s})d[X,X]_{s}^{c}\\ &\ +\sum_{0<s\le t}(f(X_{s},Y_{s})-f(X_{s-},Y_{s-})-f_{x}(X_{s-},Y_{s-})\Delta X_{s}\n
&\hskip 6cm -f_{y}(X_{s-},Y_{s-})\Delta Y_{s})\ .\nonumber
\end{align}
Since $Y$ is FV we may write
\begin{align}
\int_{(0,t]}f_y(X_{s-},Y_{s-})dY_s=\int_0^tf_y(X_s,Y_s)dY_s^c+\sum_{0<s\le t}f_y(X_{s-},Y_{s-})\Delta Y_s
\end{align}
to obtain, after cancellation of the sum part, that
\begin{align}\label{eq:ny}
f(X_{t},Y_{t})&=f(X_{0},Y_{0})+\int_{(0,t]}f_{x}(X_{s-},Y_{s-})dX_{s}+\int_{(0,t]}f_{y}(X_{s},Y_{s})dY_{s}^{c}\n
&\ +\frac{1}{2}\int_{0}^{t}f_{xx}(X_{s},Y_{s})d[X,X]_{s}^{c}\\ &\ +\sum_{0<s\le t}(f(X_{s},Y_{s})-f(X_{s-},Y_{s-})-f_{x}(X_{s-},Y_{s-})\Delta X_{s})\ .\nonumber
\end{align}
Recalling \eq{LK} and noting that
\begin{align}\label{eq:10}
\int_{\mathbb{R}\setminus[-1,1]\times(0,t]}zN(dz,ds)=\sum_{0<s\le t}\Delta X_s\one{|\Delta X_s|>1}
\end{align}
and
\begin{align}\label{eq:11}
d[X,X]^{c}_{s}&=d[\sigma W,\sigma W]_s=\sigma^{2}\, ds
\end{align}
we have that
\begin{align}\label{eq:12}
f(X_{t},Y_{t})&=f(X_{0},Y_{0})+\int_{(0,t]}\left(cf_{x}(X_{s},Y_{s})+\frac{\sigma^{2}}{2}f_{xx}(X_{s},Y_{s})\right)ds\n&\ +\int_{(0,t]}f_{y}(X_{s},Y_{s})dY_{s}^{c}\\
&\ +\sigma\int_{0}^{t}f_{x}(X_{s},Y_{s})dW_{s}+\int_{[-1,1]\times(0,t]}f_{x}(X_{s-},Y_{s-})z \tilde N(dz,ds)\n
&\ +\sum_{0<s\le t}\left(f(X_{s},Y_{s})-f(X_{s-},Y_{s-})-f_{x}(X_{s-},Y_{s-})\Delta X_{s}1_{\{|\Delta X_{s}|\le 1\}}\right)\ .\nonumber
\end{align}
Next, note that
\begin{align}
&\sum_{0<s\le t}\left(f(X_{s},Y_{s})-f(X_{s-},Y_{s-})-f_{x}(X_{s-},Y_{s-})\Delta X_{s}1_{\{|\Delta X_{s}|\le 1\}}\right)\n &=
\sum_{0<s\le t}\left(f(X_{s},Y_{s})-f(X_{s},Y_{s-})\right)\\
&\ \ +\sum_{0<s\le t}\left(f(X_{s-}+\Delta X_{s},Y_{s-})-f(X_{s-},Y_{s-})-f_{x}(X_{s-},Y_{s-})\Delta X_{s}1_{\{|\Delta X_{s}|\le 1\}}\right)\nonumber
\end{align}
and that if in addition $\int_{\mathbb{R}\setminus[-1,1]}|f(x+z,y)-f(x,y)|\nu(dz)$ is bounded (as a function of $(x,y)$) on compact sets, then
\begin{align}\label{eq:14}
&\sum_{0<s\le t}\left(f(X_{s-}+\Delta X_{s},Y_{s-})-f(X_{s-},Y_{s-})-f_{x}(X_{s-},Y_{s-})\Delta X_{s}1_{\{|\Delta X_{s}|\le 1\}}\right)\n
&=\int_{\mathbb{R}\times(0,t]}\left(f(X_{s-}+z,Y_{s-})-f(X_{s-},Y_{s-})-f_{x}(X_{s-},Y_{s-})z1_{\{|z|\le 1\}}\right)N(dz,ds)\\
&=\int_{\mathbb{R}\times(0,t]}\left(f(X_{s-}+z,Y_{s-})-f(X_{s-},Y_{s-})-f_{x}(X_{s-},Y_{s-})z1_{\{|z|\le 1\}}\right)\tilde N(dz,ds)\n
&\ \ +\int_0^t\left(\int_{\mathbb{R}}\left(f(X_{s-}+z,Y_{s-})-f(X_{s-},Y_{s-})-f_{x}(X_{s-},Y_{s-})z1_{\{|z|\le 1\}}\right)\nu(dz)\right)ds\ .\nonumber
\end{align}
If we denote
\begin{align}\label{eq:A}
\mathcal{A}f(x,y)=cf_x(x,y)+\frac{\sigma^2}{2}f_{xx}(x,y)+\int_{\mathbb{R}}
\left(f(x+z,y)-f(x,y)-f_x(x,y)z1_{\{|z|\le 1\}}\right)\nu(dz)\ ,
\end{align}
then putting everything together, noting that
\begin{align}
\int_{[-1,1]\times(0,t]}f_{x}(X_{s-},Y_{s-})z \tilde N(dz,ds)
\end{align}
in \eq{14} cancels with the corresponding term in \eq{12}, gives
\begin{align}
f(X_t,Y_t)&=f(X_0,Y_0)+\int_0^t\mathcal{A}f(X_s,Y_s)ds+\int_0^tf_y(X_s,Y_s)dY^c_s+\sum_{0<s\le t}(f(X_s,Y_s)-f(X_s,Y_{s-}))\n
&\ \ +\sigma\int_0^tf_x(X_s,Y_s)dW_s+\int_{\mathbb{R}\times(0,t]}(f(X_{s-}+z,Y_{s-})-f(X_{s-},Y_{s-}))\tilde N(dz,ds)\ .
\end{align}
Therefore,
\begin{align}\label{eq:loc}
M_t&=f(X_t,Y_t)-f(X_0,Y_0)-\int_0^t\mathcal{A}f(X_s,Y_s)ds\n &\ \ -\int_0^tf_y(X_s,Y_s)dY^c_s-\sum_{0<s\le t}(f(X_s,Y_s)-f(X_s,Y_{s-}))\\
&=\sigma\int_0^tf_x(X_s,Y_s)dW_s+\int_{\mathbb{R}\times(0,t]}(f(X_{s-}+z,Y_{s-})-f(X_{s-},Y_{s-}))\tilde N(dz,ds)\nonumber
\end{align}
is a local martingale (e.g., Subsection~4.3.2, p. 230-233 in \cite{applebaum}, Thm.~29, p. 171 of \cite{p2004} and Prop. 4.10 in \cite{r2004}).

Next, if we denote
\begin{align}
U_t&=\sigma\int_0^tf_x(X_s,Y_s)dW_s\\
V_t&=\int_{\mathbb{R}\times(0,t]}(f(X_{s-}+z,Y_{s-})-f(X_{s-},Y_{s-}))\tilde N(dz,ds)\ ,\nonumber
\end{align}
then $V$, being a compensated sum of jumps, is quadratic pure jump (e.g., \cite{meyer}) with
\begin{align}\label{eq:VV}
[V,V]_t=\int_{\mathbb{R}\times(0,t]}(f(X_{s-}+z,Y_{s-})-f(X_{s-},Y_{s-}))^2N(dz,ds)
\end{align}
(note: $N$, not $\tilde N$), $U$ has quadratic variation
\begin{align}
[U,U]_t=\sigma^2\int_0^tf_x^2(X_s,Y_s)ds
\end{align}
(e.g., Thm. 29, p. 75 of \cite{p2004}) and $[V,U]_t=0$.

Now, consider the following.
\begin{Assumption}\label{AU}
There exists a closed set $B\subset\mathbb{R}^2$ satisfying $P((X_t,Y_t)\in B)=1$ for all $t\ge 0$ such that
$f_x(x,y)$ and $\int_{\mathbb{R}}(f(x+z,y)-f(x,y))^2\nu(dz)$ are bounded on $B$.
\end{Assumption}

Under Assumption~\ref{AU}, it follows that
\begin{align}\label{eq:qv}
[M,M]_t=\int_0^t\left(\sigma^2 f_x^2(X_s,Y_s)+\int_{\mathbb{R}}(f(X_s+z,Y_s)-f(X_s,Y_s))^2\nu(dz)\right)ds +\tilde M_t
\end{align}
where
\begin{align}\label{eq:tildeM}
\tilde M_t=\int_{\mathbb{R}\times(0,t]}(f(X_{s-}+z,Y_{s-})-f(X_{s-},Y_{s-}))^2\tilde N(dz,ds)
\end{align}
is a zero mean martingale (e.g., Prop.~4.10 in \cite{r2004}, or a generalized version of Lemma~1 of \cite{ny2005}).

Therefore, under Assumption~\ref{AU}, we now have (as in \cite{kb2013}), that
\begin{align}\label{eq:L_2}
E[M,M]_t\le \int_0^t C(s)ds
\end{align}
where
\begin{align}C(s)=E\left(\sigma^2 f_x^2(X_s,Y_s)+\int_{\mathbb{R}}(f(X_s+z,Y_s)-f(X_s,Y_s))^2\nu(dz)\right)
\end{align}
is bounded (in $s$).

\begin{Lemma}\label{C}
Assume that $M$ is a local martingale for which $E[M,M]_t$ is absolutely continuous (with respect to Lebesgue measure) and has a bounded (necessarily nonnegative) density C(s). Then
\begin{description}
\item{(1)} $M$ is an $L^2$ martingale,
\item{(2)} for every (deterministic) $h$ with $h(t)\sqrt{t}\to 0$ as $t\to\infty$; $h(t)M_t\to 0$ in $L^2$, as $t\to\infty$, and
\item{(3)} for every continuous, nonnegative, nonincreasing $h$, satisfying $\int_{t_0}^\infty h^2(s)ds<\infty$ for some $t_0\ge 0$; $h(t)M_t\to0$ a.s., as $t\to\infty$.
\end{description}
In particular, for every $\gamma>1/2$, $M_t/t^\gamma\to 0$ in $L^2$ and a.s., as $t\to\infty$.
\end{Lemma}

\proof
Let $C=\sup_{s\ge 0}C(s)$. From $E[M,M]_t\le Ct<\infty$ it follows that $M$ is an $L^2$ martingale with $EM^2_t=E[M,M]_t$ (e.g., Cor.~3, p.~73 of \cite{p2004}). This implies that
\begin{equation}
E(h(t)M_t)^2=h^2(t)E[M,M]_t\le Ct\,h^2(t)\ .
\end{equation}
and thus (2) follows.

Next, we prove (3): if $h(t)=0$ for some $t>0$ then $h(s)M_s=0$ for $s\ge t$. Also, if $h(t)>1$ for some $t$ then we may replace $h$ by $h_1(t)=h(t)\wedge 1$ and clearly $h(t)M_t\to 0$ a.s. if and only if $h_1(t)M_t\to 0$ a.s. and $\int_0^\infty h_1^2(s)ds<\infty$. Thus, we may restrict ourselves to $h$ with $0<h(t)\le 1$ for every $t\ge 0$, such that $\int_0^\infty h^2(s)ds<\infty$. With this assumption,
$h\cdot M_t\equiv\int_{(0,t]}h(s)dM_s$ is a martingale with
\begin{align}\label{eq:as}
E(h\cdot M_t)^2&=E[h\cdot M,h\cdot M]_t=\int_0^th^2(s)dE[M,M]_s\n &=\int_0^th^2(s)C(s)ds\le C\int_0^\infty h^2(s)ds
\end{align}
(for the second equality see Th.~29, p.~75 of \cite{p2004}) and thus converges a.s. Consider now $A(t)=\frac{1}{h(t)}-1$. Then $A(\cdot)$ is continuous, nonnegative, nondecreasing and, as $t\to\infty$, $A(t)\to\infty$ with $\int_0^t\frac{dM_s}{1+A(s)}$ converging a.s.
Hence, as in Ex.~14, p.~95 of \cite{p2004}, we also have that $h(t)M_t=M_t/(1+A(t))\to 0$ a.s.\qed

Thus, we can now conclude the following.
\begin{Theorem}\label{mart1}
With \cadlag and adapted $X$, $Y$ and with a function $f:\mathbb{R}^2\to \mathbb{R}$, where
\begin{itemize}
\item $X$ is a real valued L\'evy process (with respect to the underlying filtration) with L\'evy triplet $(c,\sigma^2,\nu(\cdot))$,
\item $Y$ a FV process,
\item $f\in \mathcal{C}^{2,1}$ with $\int_{\mathbb{R}\setminus[-1,1]}|f(x+z,y)-f(x,y)|\nu(dz)$ bounded on compact sets (redundant under Assumption~\ref{AU})
\end{itemize}
and with $\mathcal{A}$ defined in \eq{A} then (\ref{eq:M_t})
%\begin{align}
%M_t&=f(X_t,Y_t)-f(X_0,Y_0)-\int_0^t\mathcal{A}f(X_s,Y_s)ds\n &\ \ -\int_0^tf_y(X_s,Y_s)dY^c_s-\sum_{0<s\le t}(f(X_s,Y_s)-f(X_s,Y_{s-}))
%\end{align}
is a local martingale.

If in addition Assumption~\ref{AU} holds, then the assumptions and hence the conclusions of Lemma~\ref{C} hold.
\end{Theorem}

\begin{Remark}\label{R1}\rm
We note the following regarding Assumption~\ref{AU}:
\begin{itemize}
\item A sufficient condition is that $f$ and $f_x$ are bounded on $\mathbb{R}^2$. For example, it holds for $f(x,y)=\sin(\alpha(x+y))$ or $f(x,y)=\cos(\alpha(x+y))$ and, thus, also for $f(x,y)=e^{i\alpha(x+y)}$.
\item Another sufficient condition is that $f(x+z,y)$ is bounded on $(x,y)\in B$ and $z\in\mathbb{R}\setminus[-1,1]/(0,1]/[-1,0))$ and $f_x(x+z,y)$ is bounded on $(x,y)\in B$ and $z\in [-1,1]/(0,1]/[-1,0)$ for the general/spectrally positive/spectrally negative cases, respectively. For example, it holds for the spectrally positive case where $X_t+Y_t\ge 0$, a.s., and $f(x,y)=e^{-\alpha(x+y)}$ for $\alpha>0$.

\end{itemize}
\end{Remark}

From Remark~\ref{R1}, recalling \eq{Apsi}, $\psi$ from \eq{psi}, denoting $\varphi(\alpha)=\psi(i\alpha)=\log Ee^{-\alpha X_1}$ (real valued) for the spectrally positive case and noting that
\begin{align}
f(X_{s},Y_{s})-f(X_{s},Y_{s-})=e^{i\alpha(X_s+Y_{s})}-e^{i\alpha (X_{s}+Y_{s-})}=e^{i\alpha(X_s+Y_s)}(1-e^{-i\alpha\Delta Y_s})
\end{align}
(with $-\alpha$ replacing $i\alpha$ for the spectrally positive case) we immediately reproduce the following (see \cite{kb2013,kw1992}).
\begin{Corollary}
With $Z_t=X_t+Y_t$, where $X,Y$ are as in Theorem~\ref{mart1},
\begin{align}
M_t&=\psi(\alpha)\int_0^te^{i\alpha Z_s}ds+e^{i\alpha Z_0}-e^{i\alpha Z_t}\n
&\ +i\alpha\int_0^te^{i\alpha Z_s}dY^c_s+\sum_{0<s\le t}e^{i\alpha Z_s}(1-e^{-i\alpha\Delta Y_s})
\end{align}
is an (complex valued) $L^2$ martingale that satisfies the assumptions and hence the conclusions of Lemma~\ref{C}.

If in addition $Z_t\ge 0$ a.s. for all $t\ge0$ and $\nu(-\infty,0)=0$ then the same holds for
\begin{align}
M_t&=\varphi(\alpha)\int_0^te^{-\alpha Z_s}ds+e^{-\alpha Z_0}-e^{-\alpha Z_t}\n
&\ -\alpha\int_0^te^{-\alpha Z_s}dY^c_s+\sum_{0<s\le t}e^{-\alpha Z_s}(1-e^{\alpha\Delta Y_s})\ .
\end{align}
\end{Corollary}

\section{Multivariate Jump diffusion $X$ and multivariate $Y$}\label{multijumpdiffusion}
The components of an $n$-dimensional jump-diffusion process are as follows (e.g., Chapter~6 of \cite{applebaum}, noting the footnote on p. 363).
\begin{enumerate}
\item $W=(W_1,\ldots,W_k)^T$, where $W_i$ are independent Wiener processes.
\item $N(dz,dt)$ is a Poisson random measure on $\mathbb{R}^m\times \mathbb{R}_+$ with mean measure $\nu(dz)dt$, satisfying $\int_{\mathbb{R}^m}(\|x\|^2\wedge 1)\nu(dx)<\infty$ ($\|x\|$ is the Euclidean norm) and $\nu(\{0\})=0$.
\item $\tilde N(dz,dt)=N(dz,dt)-\nu(dz)dt$.
\item $b_i:\mathbb{R}^n\to\mathbb{R}$, $\sigma_{ij}:\mathbb{R}^n\to \mathbb{R}$, $K_i:\mathbb{R}^{n+m}\to \mathbb{R}$, $i=1,\ldots,n$, $j=1,\ldots,k$. All are Borel with:
\begin{description}
\item{(i)} Lipschitz conditions: for all $x_1,x_2\in\mathbb{R}^n$, $1\le i\le n$ and $1\le j\le k$,
\begin{align*}
|b_i(x_1)-b_i(x_2)|\vee |\sigma_{ij}(x_1)-\sigma_{ij}(x_2)|\vee\left(\int_{\|z\|\le 1}(K_i(x_1,z)-K_i(x_2,z))^2\nu(dz)\right)^{1/2}\n \le \kappa\|x_1-x_2\|\ .
\end{align*}
\item{(ii)} Finiteness conditions: for each $1\le i\le n$, for some (hence all) $x\in\mathbb{R}^n$,
\begin{align}
\int_{\|z\|\le 1}K^2_i(x,z)\nu(dz)<\infty\ .
\end{align}
\item {(iii)} For each $z$ such that $\|z\|>1$ and each $1\le i\le n$, $K_i(\cdot,z)$ is continuous.
\end{description}
\end{enumerate}
$X$ is the unique (strong Markov) strong solution of:
\begin{align}\label{eq:X_i}
X_{i,t}=X_{i,0}+\int_0^tb_i(X_s)ds+\sum_{j=1}^k\int_0^t\sigma_{ij}(X_s)dW_{j,s}&+\int_{\mathbb{R}^m\times(0,t]}K_i(X_{s-},z)\one{\|z\|>1}N(dz,ds)\n
&+\int_{\mathbb{R}^m\times(0,t]}K_i(X_{s-},z)\one{\|z\|\le 1}\tilde N(dz,ds)\\
=X_{i,0}+\int_0^tc_i(X_s)ds+\sum_{j=1}^k\int_0^t\sigma_{ij}(X_s)dW_{j,s}&+\int_{\mathbb{R}^m\times(0,t]}K_i(X_{s-},z)\one{\|K(X_{s-},z)\|>1}N(dz,ds)\n
&+\int_{\mathbb{R}^m\times(0,t]}K_i(X_{s-},z)\one{\|K(X_{s-},z)\|\le 1}\tilde N(dz,ds)\nonumber
\end{align}
where in the last expression we define
\begin{align}\label{eq:cb1}
c_i(x)=b_i(x)+\int_{\|z\|\le 1}K_i(x,z)\one{\|K(x,z)\|>1}\nu(dz)-\int_{\|z\|>1}K_i(x,z)\one{\|K(x,z)\|\le 1}\nu(dz)\ ,
\end{align}
The second equality of \eq{X_i} is straightforward, although somewhat cumbersome. We note that both integrals on the right side of \eq{cb1} are finite since
\begin{align}
\int_{\|z\|\le 1}|K_i(x,z)|\one{\|K(x,z)\|>1}\nu(dz)\le \int_{\|z\|\le 1}\|K(x,z)\|\one{\|K(x,z)\|>1}\nu(dz)\le \int_{\|z\|\le 1}\|K(x,z)\|^2\nu(dz)<\infty
\end{align}
and
\begin{align}
\int_{\|z\|>1}|K_i(x,z)|\one{\|K(x,z)\|\le 1}\nu(dz)\le \nu(\{z|\|z\|>1\})<\infty\ .
\end{align}

In addition to $X$ we also introduce an $r$ dimensional FV process $Y$.

For some $f:\mathbb{R}^{n+r}\to \mathbb{R}$, by $f\in\mathcal{C}^{2,1}$ here we mean that $f$ is twice continuously differentiable in the first $n$ coordinates and continuously differentiable in the last $r$ coordinates.

For $f:\mathbb{R}^{n+r}\to \mathbb{R}$ with $f\in\mathcal{C}^{2,1}$ we note that as for \eq{ny}, \eq{10} and \eq{11} (and the same justification) we have that
\begin{align}\label{eq:nya}
f(X_{t},Y_{t})&=f(X_{0},Y_{0})+\sum_{i=1}^n\int_{(0,t]}f_{x_i}(X_{s-},Y_{s-})dX_{i,s}+\sum_{j=1}^r\int_{(0,t]}f_{y_j}(X_{s},Y_{s})dY_{j,s}^{c}\n
&\ +\frac{1}{2}\sum_{i=1}^n\sum_{j=1}^n\int_{0}^{t}f_{x_ix_j}(X_{s},Y_{s})d[X_i,X_j]_{s}^{c}\\ &\ +\sum_{0<s\le t}\left(f(X_{s},Y_{s})-f(X_{s-},Y_{s-})-\sum_{i=1}^nf_{x_i}(X_{s-},Y_{s-})\Delta X_{i,s}\right)\ ,\nonumber
\end{align}
\begin{align}\label{eq:10a}
\int_{\mathbb{R}\setminus[-1,1]\times(0,t]}K_i(X_{s-},z)\one{\|K(X_{s-},z)\|>1}N(dz,ds)=\sum_{0<s\le t}\Delta X_{i,s}\one{\|\Delta X_s\|>1}
\end{align}
and
\begin{align}\label{eq:11a}
d[X_i,X_j]^{c}_{s}&=\sum_{\ell_1=1}^k\sum_{\ell_2=1}^k\sigma_{i\ell_1}(X_s)\sigma_{j\ell_2}(X_s)d[W_{\ell_1},W_{\ell_2}]_s
=\sum_{\ell=1}^k\sigma_{i\ell}(X_s)\sigma_{j\ell}(X_s)ds\ .
\end{align}
Thus, in a similar manner as for the L\'evy case, we are led to the operator $\mathcal{A}$ in this case defined as follows, with $c=(c_i),\sigma=(\sigma_{ij}), K=(K_i)$,
\begin{align}\label{eq:gen3}
\mathcal{A}f(x)&=c(x)^T\nabla_xf(x,y)+\frac{1}{2}\text{tr}(\sigma(x)\nabla_{xx}f(x,y)\sigma(x)^T)\\ &+\int_{\mathbb{R}^m}(f(x+K(x,z),y)-f(x,y)-\nabla_x f(x,y)^TK(x,z)
\one{\|K(x,z)\|\le 1})\nu(dz)\nonumber
\end{align}
where we may replace the last integral in \eq{gen3} by
\begin{align}
\int_{\mathbb{R}^m}(f(x+z,y)-f(x,y)-\nabla_x f(x,y)^Tz
\one{\|z\|\le 1})\mu(x,dz)\ ,
\end{align}
where $\mu(x,A)=\nu(\{z|\ K(x,z)\in A\}$.

Finally consider the following.
\begin{Assumption}\label{AU2}
There exists a closed set $B\subset\mathbb{R}^{n+r}$ with $P((X_t,Y_t)\in B)=1$ for all $t\ge 0$ such that
$\|\nabla_xf(x,y)^T\sigma(x)\|$ and
\[\int_{\mathbb{R}^m}(f(x+K(x,z),y)-f(x,y))^2\nu(dz)=\int_{\mathbb{R}^m}(f(x+z,y)-f(x,y))^2\mu(x,dz)
 \]
 are bounded on $B$.
\end{Assumption}

A repetition of the arguments from Section~\ref{Levy} gives the following result.
\begin{Theorem}\label{mart3}
With c\`adl\`ag and adapted $X$, $Y$ and with a function $f:\mathbb{R}^{n+r}\to \mathbb{R}$, where
\begin{itemize}
\item $X$ is an $\mathbb{R}^n$ valued jump diffusion process (with respect to the underlying filtration) with $(b(\cdot),\sigma(\cdot),K(\cdot,\cdot),\nu(\cdot))$ satisfying (i)-(iii),
\item $Y$ a $\mathbb{R}^r$ valued FV process,
\item $f\in \mathcal{C}^{2,1}$ and
    \begin{align*}
    \int_{\mathbb{R}^m}\left|f(x+K(x,z),y)-f(x,y)\right|1_{\{\|K(x,z)\|>1\}}\nu(dz)\\
    =\int_{\mathbb{R}^m}\left|f(x+z,y)-f(x,y)\right|1_{\{\|z\|>1\}}\mu(x,dz)
    \end{align*}
    is bounded on compact sets (redundant under Assumption~\ref{AU2}),
\end{itemize}
and with $\mathcal{A}$ defined via \eq{gen3}, then \eq{M_t}
%\begin{align}
%M_t&=f(X_t,Y_t)-f(X_0,Y_0)-\int_0^t\mathcal{A}f(X_s,Y_s)ds\n &\ \ -\int_0^t\nabla_y f(X_s,Y_s)^TdY^c_s-\sum_{0<s\le t}(f(X_s,Y_s)-f(X_s,Y_{s-}))
%\end{align}
is a local martingale.

If in addition Assumption~\ref{AU2} holds, then the assumptions and hence the conclusions of Lemma~\ref{C} hold.
\end{Theorem}

\begin{Remark}\label{axes}\rm
An important special case of this setup is obtained when the measure $\nu$ is concentrated on the axes. That is, on $\bigcup_{i=1}^m\mathbf{1}_i\mathbb{R}$, where $\mathbf{1}_i$ is a unit vector with one in the $i$ coordinate and zero elsewhere, so that $\mathbf{1}_i\mathbb{R}$ denotes the $i$th axis. For this case, we let $\nu_j$ be the $j$th marginal of $\nu$ and let $K^j(x,z)=K(x,\mathbf{1}_jz)$. For this case, we have in \eq{gen3} that
\begin{align}
c_i(x)=b_i(x)+\sum_{j=1}^m\left(\int_{|z|\le 1}K_i^j(x,z)\one{\|K^j(x,z)\|>1}\nu_j(dz)-\int_{|z|>1}K_i^j(x,z)\one{\|K^j(x,z)\|\le 1}\nu_j(dz)\right)\ ,
\end{align}
the last summand becomes
\begin{align}
&\int_{\mathbb{R}^m}(f(x+K(x,z),y)-f(x,y)-\nabla_x f(x,y)^TK(x,z)
\one{\|K(x,z)\|\le 1})\nu(dz)\\
&=\sum_{j=1}^m\int_{\mathbb{R}}(f(x+K^j(x,z),y)-f(x,y)-\nabla_x f(x,y)^TK^j(x,z)
\one{\|K^j(x,z)\|\le 1})\nu_j(dz)\nonumber
\end{align}
and in Assumption~\ref{AU2} we have that
\begin{align}
\int_{\mathbb{R}^m}(f(x+K(x,z),y)-f(x,y))^2\nu(dz)=\sum_{j=1}^n\int_{\mathbb{R}}(f(x+K^j(x,z),y)-f(x,y))^2\nu_j(dz)
\end{align}
which is bounded if and only if every term on the right is bounded. We also note that in this case we can replace the first equality in \eq{X_i} by
\begin{align}
X_{i,t}=X_{i,0}+\int_0^tb_i(X_s)ds+\sum_{j=1}^k\int_0^t\sigma_{ij}(X_s)dW_{j,s}&+\sum_{j=1}^m\int_{\mathbb{R}\setminus[-1,1]\times(0,t]}K^j_i(X_{s-},z)N_j(dz,ds)\n
&+\sum_{j=1}^m\int_{[-1,1]\times(0,t]}K_i^j(X_{s-},z)\tilde N_j(dz,ds)
\end{align}
where $N_1,\ldots,N_m$ are independent Poisson random measures (on $\mathbb{R}\times\mathbb{R}_+$) with intensities $\nu_j(dz)dt$, $j=1,\ldots,m$.
\end{Remark}

This setup will prove useful for the next section.

\section{Markov additive $X$ with finite state space modulation}\label{markovadditive}
By a {\em Markov additive process with finite state space modulation} we mean a process $(J,X)$ where $J$ is a finite state space Markov chain with some rate transition matrix $Q=(q_{ij})$ and during epochs where $J(t)=i$, $X$ behaves like a L\'evy process with some triplet $(c_i,\sigma_i^2,\nu_i)$. In addition, at state change epochs of $J$ from $i$ to $j$ the process $X$ may incur independent jumps that have a distribution $G_{ij}$ when the transition is from $i$ to $j$. For a precise description see, e.g., \cite{ak2000}. Let us construct $(X,J)$ as a two dimensional jump diffusion with $n=2$, $k=1$ and $m=K(K+1)$, where the states are $1,\ldots,K$.

Let us begin with the construction of $J$. In view of Remark~\ref{axes}, for every $1\le i,j\le j$ with $i\not=j$, we let $N_{ij}$ denote independent Poisson random measures on $\mathbb{R}\times\mathbb{R}_+$ with intensities $\nu_{ij}(\cdot)=q_{ij}G_{ij}(\cdot)$. Then $J$ satisfies the following equation.
\begin{align}
J_t=J_0+\sum_{i\not=j}\int_{(0,t]}\int_\mathbb{R} (j-i)\one{J_{s-}=i}N_{ij}(dz,ds)
\end{align}
Clearly one can find functions $K_1^{ij}((\cdot,\cdot),\cdot)$ that satisfy (i)-(iii) such that $K_1^{ij}((u,x),z)=(j-i)\one{u=i}$ for all $i,j,u,x,z$.

 Now let $N_i$ being independent Poisson random measures (also on $\mathbb{R}\times\mathbb{R}_+$). We also assume that they are independent of $\{N_{ij}|1\le i,j\le K\}$ with intensities $\nu_i(dz)dt$, where $\nu_i$ is a L\'evy measure and $W$ is an independent Wiener process. Then $X$ satisfies
\begin{align}
X_t=X_0+\int_0^tc_{J_s}ds+\int_0^t\sigma_{J_s}dW_s&+\sum_{i=1}^K \int_{(0,t]}\int_{\mathbb{R}\setminus[-1,1]}\one{J_{s-}=i}z N_i(dz,ds)\n
&+\sum_{i=1}^K \int_{(0,t]}\int_{[-1,1]}\one{J_{s-}=i}z \tilde N_i(dz,ds)
\\
&+\sum_{i\not=j} \int_{(0,t]}\int_{\mathbb{R}}\one{J_{s-}=i}zN_{ij}(dz,ds)\nonumber
\end{align}
where here also one can find $K^{ij}_2((\cdot,\cdot),\cdot)$, $K_2^i((\cdot,\cdot),\cdot)$,$c_2(\cdot,\cdot)$ and $\sigma_2(\cdot,\cdot)$ which satisfy (i)-(iii) and that agree with the corresponding values in the equation.

Noting that it is easy to construct $g(z,x,y)\in\mathcal{C}^{2,1}$ such that $g(i,x,y)=f(i,x,y)$ and, for each $(x,y)$, $g(\cdot,x,y)$ is constant in a neighborhood of $i$, then for appropriate $f$ (see Theorem~\ref{mart4}), the operator $\mathcal{A}$ in this case becomes
\begin{align}
\mathcal{A}f(i,x,y)=c_i f_x(i,x,y)&+\frac{1}{2}\sigma_i^2f_{xx}(i,x,y)\n &+\int_{\mathbb{R}}(f(i,x+z,y)-f(i,x,y)-f_x(i,x,y)z\one{|z|\le 1})\nu_i(dz)\\
&+\sum_{j\not=i}q_{ij}\int_{\mathbb{R}}(f(j,x+z,y)-f(i,x,y))G_{ij}(dz)\ ,\nonumber
\end{align}
noting that for this case there is no need to start with $b_i(x)$ and then modify to $c_i(x)$ as was done in~\eq{cb1}.
If we recall that $q_{ii}=-\sum_{j\not=i}q_{ij}$ and denote $G_{ii}(\{0\})=1$, then we may rewrite $\mathcal{A}$ as follows
\begin{align}\label{eq:gen4}
\mathcal{A}f(i,x,y)=c_i f_x(i,x,y)&+\frac{1}{2}\sigma_i^2f_{xx}(i,x,y)\n &+\int_{\mathbb{R}}(f(i,x+z,y)-f(i,x,y)-f_x(i,x,y)z\one{|z|\le 1})\nu_i(dz)\\
&+\sum_{j=1}^Kq_{ij}\int_{\mathbb{R}}f(j,x+z,y)G_{ij}(dz)\ .\nonumber
\end{align}

Now consider
\begin{Assumption}\label{AU3}
There exists a closed set $B\subset\{1,\ldots,K\}\times\mathbb{R}^2$ with $P((J_t,X_t,Y_t)\in B)=1$ for all $t\ge 0$ such that for each $1\le i\le K$,
$f_x(i,x,y)$ and
$\int_{\mathbb{R}}(f(i,x+z,y)-f(i,x,y))^2\nu_i(dz)$
 are bounded on $B$.
\end{Assumption}

Together with Theorem~\ref{mart3} and observing that no continuity or differentiability with respect to the first variable is needed we now have the following.

\begin{Theorem}\label{mart4}
With c\`adl\`ag and adapted $J$, $X$, $Y$ and with a function $f:\mathbb{R}^3\to \mathbb{R}$, where
\begin{itemize}
\item $(J,X)$ is a Markov additive process as described above,
\item $Y$ an $\mathbb{R}^r$-valued FV process,
\item $f(i,\cdot,\cdot)\in \mathcal{C}^{2,1}$ and $\int_{\mathbb{R}\setminus[-1,1]}|f(i,x+z,y)-f(i,x,y)|\nu(dz)$ is bounded on compact intervals (redundant under Assumption~\ref{AU3}) for each $1\le i\le K$
\end{itemize}
and with $\mathcal{A}$ defined via \eq{gen4}, then \eq{M_t}, with $(J,X)$ replacing $X$,
%\begin{align}
%M_t&=f(J_t,X_t,Y_t)-f(J_0,X_0,Y_0)-\int_0^t\mathcal{A}f(J_s,X_s,Y_s)ds\n &\ \ -\int_0^t \nabla_y f(J_s,X_s,Y_s)^TdY^c_s-\sum_{0<s\le t}(f(J_s,X_s,Y_s)-f(J_s,X_s,Y_{s-}))
%\end{align}
is a local martingale.

If in addition Assumption~\ref{AU3} holds, then the assumptions and hence the conclusions of Lemma~\ref{C} hold.
\end{Theorem}

We note that if we take $f(i,x,y)=g(x,y)\delta_{ii_0}$ for some $i_0$ ($\delta_{ii_0}=1$ for $i=i_0$ and zero otherwise), and
\begin{align}
\mathcal{A}_jg(x,y)=c_jg_x(x,y)+\frac{\sigma_j^2}{2}g_{xx}(x,y)+\int_{\mathcal{R}}
\left(g(x+z,y)-g(x,y)-g_x(x,y)z1_{\{|z|\le 1\}}\right)\nu_j(dz)
\end{align}
then $\mathcal{A}$ becomes
\begin{align}\label{eq:gen5}
\mathcal{A}f(i,x,y)=\mathcal{A}_{i_0}g(x,y)\delta_{ii_0}+q_{ii_0}
\int_{\mathbb{R}}g(x+z,y)G_{ii_0}(dz)
\end{align}
or if we prefer to write this in matrix notation where for each $ij$ we compute the operator of $g(x,y)\one{i=j}$ then we have the following matrix valued operator $\mathcal{F}$ given by
\begin{equation}
\mathcal{F}g(x,y)=\diag{\mathcal{A}_1g(x,y),\ldots,\mathcal{A}_Kg(x,y)}+Q\circ \int_{\mathbb{R}}g(x+z,y)G(dz)
\end{equation}
where $A\circ B=(a_{ij}b_{ij})$ and $\int_{\mathbb{R}}g(x+z,y)G(dz)=\left(\int_{\mathbb{R}}g(x+z,y)G_{ij}(dz)\right)$. Finally, recalling the notation $\mathbf{1}_i$, we now have the following.

\begin{Corollary}
 With the assumptions of Theorem~\ref{mart3}, if $g(x,y)\in \mathcal{C}^{2,1}$ where $\int_{\mathbb{R}\setminus[-1,1]}|g(x+z,y)-g(x,y)|\nu(dz)$ is bounded on compact sets (redundant under Assumption~\ref{AU3}), then the following is a $K$-dimensional local martingale (a vector of local martingales)
\begin{align}
\mathbf{M}_t=g(X_t,Y_t)\mathbf{1}^T_{J_t}-g(X_0,Y_0)\mathbf{1}^T_{J_0}-\int_0^t\mathbf{1}_{J_s}^T\mathcal{F}g(X_s,Y_s)ds\n
-\sum_{i=1}^r\int_0^tg_{y_i}(X_s,Y_s)\mathbf{1}_{J_s}^TdY^c_{i,s}+\sum_{0<s\le t}(g(X_{s},Y_{s})-g(X_s,Y_{s-}))\mathbf{1}_{J_s}^T\ .
\end{align}

If in addition $g_x(x,y)$ and $\int_{\mathbb{R}}(g(x+z,y)-g(x,y))^2\nu_i(dz)$ are bounded on $B$ from Assumption~\ref{AU3} for each $1\le i\le K$, then the assumptions and hence the conclusions of Lemma~\ref{C} hold for each coordinate of $\mathbf{M}$.
\end{Corollary}

We note that a special case of this last result was introduced in \cite{ak2000} for the case where $Y$ is one dimensional and continuous, $g(x,y)=e^{\alpha(x+y)}$ under various restrictions on $\alpha$ (depending on whether the L\'evy processes involved are general, spectrally positive with nonnegative jumps at state change epochs or spectrally negative with nonpositive jumps at state chage epochs). In this case it is easy to check that
\begin{align}
\mathcal{F}g(x,y)= e^{\alpha(x+y)}F(\alpha)
\end{align}
where
\begin{equation}
F(\alpha)=\diag{\psi_1(\alpha),\ldots,\psi_K(\alpha)}+Q\circ\int_{\mathbb{R}}e^{\alpha z}G(dz)
\end{equation}
and
\begin{equation}
\psi_i(\alpha)=c_i\alpha+\frac{\sigma_i^2}{2}\alpha^2+\int_{\mathbb{R}}(e^{\alpha z}-1-\alpha z\one{|z|\le 1}\nu_i(dz)
\end{equation}
 are the L\'evy exponents.

The above substantially generalizes the results in \cite{ak2000}  and in particular we have the following.
\begin{Corollary}
Under the assumptions of Theorem~\ref{mart3}, with $r=1$, the notations above and $Z_t=X_t+Y_t$, if either
\begin{enumerate}
\item $\mathfrak{R}(\alpha)=0$,
\item $Z\ge 0$ a.s., $\nu_i(-\infty,0)=0$ for all $1\le i\le K$ and $\mathfrak{R}(\alpha)\le 0$,
\item $Z\le 0$ a.s., $\nu_i(0,\infty)=0$ for all $1\le i\le K$ and $\mathfrak{R}(\alpha)\ge 0$,
\end{enumerate}
Then,
\begin{align}
\mathbf{M}_t&=e^{\alpha Z_t}\mathbf{1}_{J_t}^T-e^{\alpha Z_0}\mathbf{1}_{J_0}^T-\int_0^te^{\alpha Z_s}\mathbf{1}_{J_s}^TdsF(\alpha)\n
&\ \ \ \ -\alpha\int_0^te^{\alpha Z_s}\mathbf{1}_{J_s}^TdY^c_s-\sum_{0<s\le t}e^{\alpha Z_s}(1-e^{-\alpha \Delta Y_s})\mathbf{1}_{J_s}^T
\end{align}
is a zero mean vector valued $L^2$ martingale, of which each coordinate satisfies the assumptions and hence the conclusions of Lemma~\ref{C}.
\end{Corollary}

%\section{A final remark}\label{finalremark}
%Throughout the paper we have seen that under some conditions which in applications are typically relatively easy to verify, we have that $M_t$ is %an $L^2$ martingale satisfying the assumptions and hence the conclusions of Lemma~\ref{C}. The question is what happens for $h(t)=1/\sqrt{t}$. In %general not much can be said, but we remark that under some considerably more restrictive conditions it is plausible that $M_t/\sqrt{t}$ may %converge in distribution to a normal random variable (e.g., \cite{vz2000}).

\end{document}